\documentclass[12pt]{article}

\usepackage[a4paper,margin=3cm,innermargin=3cm]{geometry}
\usepackage{marginnote, etoolbox, imakeidx, needspace}
\usepackage{hyperref}
\usepackage[T1]{fontenc}
\usepackage[english]{babel}
\usepackage{indentfirst}

\usepackage{verbatim}

\usepackage{xcolor}

\usepackage{amsmath, amssymb, amsfonts}

\newtheorem{theorem}{Theorem}
\newtheorem{lemma}{Lemma}

\newcommand{\proof}{\textit{Proof.}}
\newcommand{\qed}{$\Box$}

\newcommand{\op}{\operatorname}
\newcommand{\abs}[1]{\left\lvert #1 \right\rvert}

\newcommand{\br}[1]{\left\{ #1 \right\}}
\newcommand{\setdef}[2]{\br{#1 \colon #2}} 
\newcommand{\pr}[1]{\left( #1 \right)} 
\newcommand{\floor}[1]{\left\lfloor #1 \right\rfloor}
\newcommand{\ceil}[1]{\left\lceil #1 \right\rceil}
\newcommand{\irange}[2]{\left[#1,#2\right]} 
\newcommand{\Pb}[1]{\mathbb{P}{\pr{#1}}}
\newcommand{\E}[1]{\op{E}\!\pr{#1}}
\newcommand{\D}[1]{\op{Var}\!\pr{#1}}

\allowdisplaybreaks

\title{The maximum size of an induced forest \\ in the binomial random graph}

\author{Akhmejanova Margarita \footnote{Computer, Electrical and Mathematical Sciences and Engineering Division, King Abdullah University of Science and Technology (KAUST), Thuwal 23955-6900, Saudi Arabia, margarita.akhmejanova@kaust.edu.sa.}, Vladislav Kozhevnikov \footnote{Moscow Institute of Physics and Technology, Laboratory of Combinatorial and Geometric Structures, Dolgoprudny, Russia, vladislavkozhevnikov@gmail.com.}}

\date{}

\begin{document}

\maketitle


\begin{abstract}
The  celebrated Frieze's result about the independence number of $G(n,p)$ states that it is concentrated in an interval of size $o(1/p)$ for all $C_{\varepsilon}/n<p=o(1)$. We show concentration
in an interval of size $o(1/p)$ for the maximum size (number of vertices) of an induced forest in $G(n,p)$ for all $C_{\varepsilon}/n<p<1-\varepsilon$. Presumably, it is the first generalization of Frieze's result to another class of induced subgraphs for such a range of~$p$.
\end{abstract}

\section{Introduction}

For a graph $G$, we denote the maximum number of vertices in its induced forest and tree by $F(G)$ and $T(G)$, respectively. Here, we prove the following result.

\begin{theorem}\label{th:concentration_max_ind_rooted_forests}
Let $p=p(n)\in(0,1)$, $q=1/(1-p)$.
Then for any fixed $\varepsilon>0$, there exists a sufficiently large constant $C_{\varepsilon}>0$ such that if $C_{\varepsilon}/n<p<1-\varepsilon$, then a.a.s. (asymptotically almost surely)
\begin{equation*}
\floor{2\log_q(enp(1-\varepsilon))+3}\le
F(G(n,p))
\le \ceil{2\log_q(enp(1+\varepsilon))+3},
\end{equation*}
i.e. $F(G(n,p))$ is concentrated in the interval of size $o(1/p)$.
\end{theorem}


The theorem is similar to the analogous theorem for independent sets given in \cite{Frieze1990}. It refines the estimate of the deviation of the size of the largest induced forest from its asymptotic value, which is already known \cite{Glock2021}. The larger the value of $p$, the better estimate is given by the theorem. In the extreme case of $p=\op{const}$ it gives the concentration of the size in three points. The proof leverages a combination of the second moment method and Talagrand's inequality.

The study of induced trees in $G(n, p)$ was initiated by Erd\H{o}s and Palka \cite{ERDOS1983145}. They showed that if $p=\op{const}$, then for  every $\varepsilon>0$ a.a.s.
 $$(2-\varepsilon)\log_{q}np\leq T(G(n,p))\leq (2+\varepsilon)\log_{q}np.$$
Fernandez de la Vega \cite{delaVega1986} considered the case where the average degree is constant. He proved that a.a.s. such graphs ($G(n,p=\frac{C}{n}),C=\op{const}$) have induced trees of size linear in $n$. This result was also independently verified by Frieze and Jackson \cite{FRIEZE1987181}, Ku{\v{c}}era and R\"{o}dl \cite{article}, and \L{}uczak and Palka \cite{LUCZAK1988257}. Later Fernandez de la Vega revisited the case $G\pr{n,p=\frac{C}{n}}$  in \cite{delaVega1996} and proved, in his own words, "nearly best possible" result for large constant $C$ and any fixed $\varepsilon>0$:
$$\frac{2 n(\ln C-\ln \ln C-1)}{C}\le T(G(n,p))\le \frac{(1+\varepsilon)2n\ln C}{C}.$$
In \cite{Palka}, Palka and Ruci{\'{n}}ski considered the case $p=\frac{C \ln n }{n}, C>e$, for which they established that for any fixed $\varepsilon>0$ a.a.s.
$$\left(1/C-\varepsilon\right)\frac{n(\ln \ln n)}{\ln n} \leq T(G(n,p)) \leq\left(2/C+\varepsilon\right)\frac{n(\ln \ln n)}{\ln n}.$$
In 2018, Dutta and Subramanian \cite{Dutta2018} for $p \geq n^{-1 / 2}(\ln n)^2$ proved that a.a.s 
$$T(G)=2\left(\log _q n p\right)+O(1 / \ln q).$$
It is worth mentioning that the range of $p$ is crucial. For example, for $n^{-2/3+\varepsilon}<p<1/\ln^2 n$  an independence number of $G(n,p)$ has already become concentrated on two values, as seen in the recent result of Bohman and Hofstad \cite{bohman2022twopoint}.
We also mention that for the case $p=\op{const}$, 2-point concentration results for $T(G)$ and $F(G)$ were achieved by Kamaldinov, Skorkin, Zhukovskii \cite{Kamaldinov2021} and
Krivoshapko, Zhukovskii \cite{Krivoshapko2021}, respectively. Also recently, Dragani{\'{c}}, Glock 
and Krivelevich \cite{Dragani2022} proved that for any $\varepsilon>0$, there exists $C_{\varepsilon}>0$ such that a.a.s. $G(n, p)$ contains an induced \textit{linear forest} ( a forest with a maximum degree of at most $2$) of size at least $(2-\varepsilon)\log_q(np)$ and components of size $(np)^{1/2}/(\ln np)^4$ whenever $C_{\varepsilon}/n\le p \le n^{-1/2}(\ln n)^2$.

\section{Induced forests in random graphs}

Let $p=p(n)\in(0,1)$ be arbitrary, $q:=1/(1-p)$.
Let $X_k$ be the number of induced rooted subforests in $G(n,p)$ on $k$ vertices.
Then $F(G(n,p))=\max\setdef{k\in\irange{1}{n}}{X_k>0}$.
\begin{lemma}
\label{lemma:expectation_n_ind_rooted_forests}
\begin{equation}
\label{eq:expectation_n_ind_rooted_forests}
\begin{aligned}
\E{X_k}
&=\binom{n}{k}q^{-\binom{k}{2}}(kpq+1)^{k-1}\\
&=\pr{enpq^{-\frac{k-3}{2}}\exp\pr{O\pr{\frac{k/n}{(1-k/n)^2}}+O\pr{\frac{\ln{k}}{k}}+O\pr{\frac{1}{kp}}}}^{k}.
\end{aligned}
\end{equation}
\end{lemma}
\proof{}
Using the expression for the number of labelled rooted forests (see \cite{Moon1970}, formula (3.4)), we obtain
$$
\begin{aligned}
\E{X_k}
&=\binom{n}{k}\sum\limits_{m=1}^{k}
\binom{k-1}{m-1}k^{k-m}p^{k-m}(1-p)^{\binom{k}{2}-(k-m)}\\
&=\binom{n}{k}q^{-\binom{k}{2}}
\sum\limits_{m=1}^{k}\binom{k-1}{m-1}(kpq)^{k-m}\\
&=\binom{n}{k}q^{-\binom{k}{2}}(kpq+1)^{k-1}.\\
\end{aligned}
$$
Applying Stirling's formula, we obtain

\begin{equation}
\label{Stirling}
\begin{aligned}
\binom{n}{k}=\pr{\frac{en}{k}\exp\pr{O\pr{\frac{k/n}{(1-k/n)^2}}+O\pr{\frac{\ln{k}}{k}}}}^{k}\\
\end{aligned}
\end{equation}

and
$$
\begin{aligned}
\E{X_k}&=\pr{\frac{en}{k}\exp\pr{O\pr{\frac{k/n}{(1-k/n)^2}}+O\pr{\frac{\ln{k}}{k}}}q^{-\frac{k-1}{2}}(kpq+1)^{1-1/k}}^{k}\\
&=\pr{enpq^{-\frac{k-3}{2}}\exp\pr{O\pr{\frac{k/n}{(1-k/n)^2}}+O\pr{\frac{2\ln{k}}{k}}+O\pr{\frac{1}{kp}}}}^{k}.
\end{aligned}
$$

\qed{}

Now, for constant $\varepsilon>0$, let
$$
\begin{aligned}
k_{+\varepsilon}&:=\ceil{2\log_q(enp(1+\varepsilon))+3},\\
k_{-\varepsilon}&:=\floor{2\log_q(enp(1-\varepsilon))+3}.
\end{aligned}
$$

\begin{lemma}\label{lemma:expectation_n_ind_forests_at_k_pm_eps}
For any fixed $\varepsilon>0$ there exists a sufficiently large constant $C_{\varepsilon}>0$ such that if $C_{\varepsilon}/n<p<1-\varepsilon$, then
$\E{X_{k_{+\varepsilon}}}\to0$ and $\E{X_{k_{-\varepsilon}}}\to+\infty$.
\end{lemma}
\proof{}
It is easy to see that if $C_{\varepsilon}/n<p<1-\varepsilon$ and $C_{\varepsilon}$ is sufficiently large, then
$$
\begin{gathered}
\frac{k_{\pm\varepsilon}}{n}
\le\frac{2\ln(np)+7}{np}\le\frac{3\ln{C_{\varepsilon}}}{C_{\varepsilon}},\\
k_{\pm\varepsilon}p\ge2(1-p)\ln(np(1-\varepsilon))
>2\varepsilon\ln(C_{\varepsilon}(1-\varepsilon))
>\varepsilon\ln{C_{\varepsilon}},\\
k_{\pm\varepsilon}\to+\infty,
\end{gathered}
$$
and therefore (again, if $C_{\varepsilon}$ is sufficiently large)
$$
\abs{\exp\pr{O\pr{\frac{k_{\pm\varepsilon}/n}{(1-k_{\pm\varepsilon}/n)^2}}+O\pr{\frac{\ln{k_{\pm\varepsilon}}}{k_{\pm\varepsilon}}}+O\pr{\frac{1}{k_{\pm\varepsilon}p}}}-1}\le\frac{\varepsilon}{2}
$$
and
$$
\begin{aligned}
\E{X_{k_{+\varepsilon}}}&\le
\pr{enpq^{-\frac{k_{+\varepsilon}-3}{2}}(1+\varepsilon/2)}^{k_{+\varepsilon}}\le
\pr{\frac{1+\varepsilon/2}{1+\varepsilon}}^{k_{+\varepsilon}}
\to0,\\
\E{X_{k_{-\varepsilon}}}&\ge
\pr{enpq^{-\frac{k_{-\varepsilon}-3}{2}}(1-\varepsilon/2)}^{k_{-\varepsilon}}\ge
\pr{\frac{1-\varepsilon/2}{1-\varepsilon}}^{k_{-\varepsilon}}
\to+\infty.
\end{aligned}
$$
\qed{}

For
$x\ge0$, $y>0$, $z>0$, $\alpha>0$, let
\begin{equation}\label{g}
g(x,y,z,\alpha):=\pr{\frac{yz^{x}}{x^{\alpha}}}^{x}.\end{equation}

\bigskip
Here we give an auxiliary lemma, which we  will use later.

\begin{lemma}\label{Lemma_g}
For any $B\ge{A}\ge0$,
\begin{equation}
\label{eq:function_g_upper_bound_0}
\max\limits_{A\le{x}\le{B}}g(x,y,z,\alpha)
\le\max\br{\exp\pr{\frac{\alpha}{2}y^{\frac{1}{\alpha}}},g(A,y,z,\alpha),g(B,y,z,\alpha)}.
\end{equation}
If, additionally, $A\ge{y^{\frac{1}{\alpha}}}$, then
\begin{equation}
\label{eq:function_g_upper_bound_1}
\max\limits_{A\le{x}\le{B}}g(x,y,z,\alpha)
\le\max\br{\pr{y\pr{\frac{e}{A}}^{\alpha}}^{\frac{A}{2}},g(A,y,z,\alpha),g(B,y,z,\alpha)}.
\end{equation}
\end{lemma}
\proof{}
Let's fix $y,z,\alpha$ and let $g(x):=g(x,y,z,\alpha)$. Then
$$
\frac{\partial}{\partial{x}}\ln{g(x)}=\frac{\partial}{\partial{x}} x (-\alpha\ln{x} + \ln{y} + x \ln{z})=\ln\pr{\frac{yz^{2x}}{(ex)^{\alpha}}}.
$$
Consider the equation
$$
\frac{\partial}{\partial{x}}\ln{g(x)}=0.
$$
If it has no solutions for $A\le{x}\le{B}$, then
$$
\max\limits_{A\le{x}\le{B}}g(x)
=\max\br{g(A),g(B)},
$$
which implies both \eqref{eq:function_g_upper_bound_0} and \eqref{eq:function_g_upper_bound_1}.

Now, assume that it has at least one solution for $A\le{x}\le{B}$ and denote the set of its solutions by $\mathcal{C}$. Then, for each $x\in\mathcal{C}$:
$$
\ln\pr{\frac{yz^{2x}}{(ex)^{\alpha}}}=0.
$$
We can rewrite it as
$$
\pr{\frac{yz^{x}}{x^{\alpha}}}^2=y\pr{\frac{e}{x}}^{\alpha}.
$$
Thus, if determine a new function $h(x)$ as 
$$
h(x):=\pr{y\pr{\frac{e}{x}}^{\alpha}}^{\frac{x}{2}},
$$
we get that $g(x)=h(x)$ for all $x\in\mathcal{C}$. The derivative
$$
\frac{\partial}{\partial{x}}\ln{h(x)}=
\frac{1}{2}\ln\pr{\frac{y}{x^{\alpha}}}
$$
has exactly one positive zero $x_m=y^{\frac{1}{\alpha}}$, which is the maximum of $h(x)$ for $x\ge0$. Thus,
$$
\max\limits_{x\in\mathcal{C}}g(x)
\le\max\limits_{x\ge0}h(x)=h(x_m)
=\exp\pr{\frac{\alpha}{2}y^{\frac{1}{\alpha}}},
$$
which implies \eqref{eq:function_g_upper_bound_0}. If we additionally assume that $A\ge{x_m}$, then
$$
\max\limits_{x\in\mathcal{C}}g(x)
\le\max\limits_{x\ge{A}}h(x)=h(A),
$$
which implies \eqref{eq:function_g_upper_bound_1}.
\qed{}

\begin{lemma}
\label{lemma:variance_n_ind_forests_at_k_m_eps}
For any fixed $\varepsilon>0$ there exists a sufficiently large constant $C_{\varepsilon}>0$ such that if $C_{\varepsilon}/n<p<1-\varepsilon$, then
\begin{equation*}
\frac{\D{X_{k_{-\varepsilon}}}}{(\E{X_{k_{-\varepsilon}}})^2}
\le\exp\pr{O\pr{\frac{(\ln{np})^4}{np}k_{-\varepsilon}}}.
\end{equation*}
If, additionally, $p\ge(\ln{n})^2/\sqrt{n}$, then
\begin{equation*}
\frac{\D{X_{k_{-\varepsilon}}}}{(\E{X_{k_{-\varepsilon}}})^2}
=o(1).
\end{equation*}
\end{lemma}
\proof{}
Let
$k=k_{-\varepsilon}$ and
$C/n<p<1-\varepsilon$, where $C$ is some constant, whose value will be tacitly assumed to be sufficiently large whenever necessary throughout the proof. Then the value $C_{\varepsilon}$ can be chosen to be the maximum assumed value of $C$. In particular, $C$ should be large enough for
the claim of Lemma~\ref{lemma:expectation_n_ind_forests_at_k_pm_eps} to hold and for $k$ to satisfy the inequalities
$$
\varepsilon\ln{np}\le
kp
\le3\ln{np}.
$$
Now, for $\D{X_k}$,
$$
\D{X_k}\le\sum\limits_{\ell=1}^{k}{F_{\ell}},
$$
where $F_{\ell}$ is the expected value of the number of (ordered) pairs
of induced rooted forests of size $k$ intersecting in $\ell$ vertices. Each such pair of forests can obtained as follows:
\begin{itemize}
\item Choose three subsets of vertices: for given two rooted forests $F_1, F_2$ sharing $l$ vertices and their common $l$ vertices,
\item choose an induced subgraph  $F=\{T_1,\ldots,T_m\}$, where each $T_i$ is a tree size of $f_i,$ 
\item  choose a rooted forest $F_1$ on $k$ vertices with $h_1$  components, containing the induced subgraph $F$. The number of such forests is represented by $f(k,h_1,\br{f_1,\ldots,f_m})$ derived in Appendex \ref{sub_forests_with_indiced_subforest}.
\item choose a rooted forest $F_2$ on $k$ vertices with $h_2$  components, containing the induced subgraph $F$. 
\end{itemize}

\begin{align}\label{main_form}
\begin{split}
F_{\ell}=\,
&\binom{n}{\ell,k-\ell,k-\ell,n-2k+\ell}
\sum\limits_{m=1}^{\ell}\frac{1}{m!}
\sum\limits_{f_1,\ldots,f_m=1}^{+\infty}\binom{\ell}{f_1,\ldots,f_m}
f_1^{f_1-2}\cdot\ldots\cdot{f_m^{f_m-2}}\\
&\sum\limits_{h_1,h_2=1}^{k-\ell+m}f(k,h_1,\br{f_1,\ldots,f_m})\cdot{f(k,h_2,\br{f_1,\ldots,f_m})}\\
&p^{2k-h_1-h_2-\ell+m}(1-p)^{2\binom{k}{2}-\binom{\ell}{2}-(2k-h_1-h_2-\ell+m)}.
\end{split}
\end{align}

Let us make auxiliary calculations:
\begin{align*}\label{aux_form}
\sum\limits_{h=1}^{k-\ell+m}&f(k,h,\br{f_1,\ldots,f_m})(pq)^{-h}=\\
&=f_1\cdot\ldots\cdot{f_m}
\sum\limits_{k_0=0}^{k-\ell}
\binom{k-\ell}{k_0}\ell^{k-\ell-k_0}
\sum\limits_{h=1}^{k_0+m}\binom{k_{0}+m-1}{h-1}(k-\ell)^{k_{0}+m-h}(pq)^{-h}
\end{align*}

Using that $\sum_{k=1}^{n} \binom{n-1}{k-1}a^{n-k}b^k=b(a+b)^{n-1}$, we get

\begin{align}
\begin{split}
&=f_1\cdot\ldots\cdot{f_m}
\sum\limits_{k_0=0}^{k-\ell}
\binom{k-\ell}{k_0}\ell^{k-\ell-k_0}
\pr{k-\ell+\frac{1}{pq}}^{k_{0}+m-1}\frac{1}{pq}\\
&={f_1}\cdot\ldots\cdot{f_m}
\pr{k+\frac{1}{pq}}^{k-\ell}
\pr{k-\ell+\frac{1}{pq}}^{m-1}\frac{1}{pq},
\end{split}
\end{align}

\begin{align}\label{aux_form_2}
\begin{split}
\sum\limits_{f_1,\ldots,f_m=1}^{+\infty}\binom{\ell}{f_1,\ldots,f_m}f_1^{f_1}\cdot\ldots\cdot{f_m^{f_m}}
&\le\sum\limits_{f_1,\ldots,f_m=1}^{+\infty}\binom{\ell}{f_1,\ldots,f_m}f_1^{f_1-1}\cdot\ldots\cdot{f_m^{f_m-1}}
\pr{\frac{\ell}{m}}^{m}\\
&=m!\binom{\ell-1}{m-1}\ell^{\ell-m}\pr{\frac{\ell}{m}}^{m}
\le\frac{\ell^{\ell+m}}{m!},
\end{split}
\end{align}

Now we put (\ref{aux_form}) and (\ref{aux_form_2}) into (\ref{main_form}),
\begin{align*}
&\begin{aligned}
F_\ell=\,&\binom{n}{\ell,k-\ell,k-\ell,n-2k+\ell}
\sum\limits_{m=1}^{\ell}\frac{1}{m!}
\sum\limits_{f_1,\ldots,f_m=1}^{+\infty}\binom{\ell}{f_1,\ldots,f_m}
f_1^{f_1}\cdot\ldots\cdot{f_m^{f_m}}\\
&\pr{k+\frac{1}{pq}}^{2(k-\ell)}\pr{k-\ell+\frac{1}{pq}}^{2(m-1)}
(pq)^{2k-\ell+m-2}q^{\binom{\ell}{2}-2\binom{k}{2}}\le
\end{aligned}\\
&\begin{aligned}
\phantom{F_{\ell}}\,&\binom{n}{\ell,k-\ell,k-\ell,n-2k+\ell}
\sum\limits_{m=1}^{\ell}
\frac{l^{l+m}}{(m!)^2}\\
&\pr{k+\frac{1}{pq}}^{2(k-\ell)}\pr{k-\ell+\frac{1}{pq}}^{2(m-1)}
(pq)^{2k-\ell+m-2}q^{\binom{\ell}{2}-2\binom{k}{2}}\\
\phantom{F_{\ell}}=\,&\binom{n}{\ell,k-\ell,k-\ell,n-2k+\ell}\pr{k+\frac{1}{pq}}^{2(k-\ell)}
\ell^{\ell}(pq)^{2k-\ell-2}q^{\binom{\ell}{2}-2\binom{k}{2}}\\
&\sum\limits_{m=1}^{\ell}
\frac{(\ell pq)^{m}}{(m!)^2}\pr{k-\ell+\frac{1}{pq}}^{2(m-1)}
\end{aligned}\\
&\begin{aligned}
\phantom{F_{\ell}}=\,&\binom{n}{\ell,k-\ell,k-\ell,n-2k+\ell}
\pr{k+\frac{1}{pq}}^{2(k-\ell)}\ell^{\ell}{(pq)^{2k-\ell-2}}
q^{\binom{\ell}{2}-2\binom{k}{2}}\\
&\sum\limits_{m=0}^{\ell-1}\frac{(\ell{pq})^{m+1}}{((m+1)!)^2}
\pr{k-\ell+\frac{1}{pq}}^{2m}
.
\end{aligned}
\end{align*}

Let's consider the following cases:
\begin{enumerate}
    \item ${k}-\ell\le\ceil{\frac{k}{\sqrt{kpq}}-\frac{1}{pq}}$;
    \item ${k}-\ell\ge\ceil{\frac{k}{\sqrt{kpq}}-\frac{1}{pq}}$.
\end{enumerate}

1. In this case
$$
\sum\limits_{m=0}^{\ell-1}\frac{1}{((m+1)!)^2}
\pr{k-\ell+\frac{1}{pq}}^{2m}
(\ell{pq})^{m+1}\le
\exp\pr{2\pr{k-\ell+\frac{1}{pq}}\sqrt{kpq}},
$$
$$
\begin{aligned}
F_{\ell}&\le\binom{n}{\ell,k-\ell,k-\ell,n-2k+\ell}
\pr{k+\frac{1}{pq}}^{2(k-\ell)}\ell^{\ell}{(pq)^{2k-\ell-2}}
q^{\binom{\ell}{2}-2\binom{k}{2}}e^{2\pr{k-\ell+\frac{1}{pq}}\sqrt{kpq}}\\
&\leq\binom{n}{\ell,k-\ell,k-\ell,n-2k+\ell}
(kpq+1)^{2k-\ell}(pq)^{-2}
q^{\binom{\ell}{2}-2\binom{k}{2}}
e^{2\pr{k-\ell+\frac{1}{pq}}\sqrt{kpq}}\\
&=\E{Y_k}\cdot\binom{k}{k-\ell}\binom{n-k}{k-\ell}(kpq+1)^{k-\ell+1}q^{-\frac{(k-\ell)(k+\ell-1)}{2}}
e^{2(k-\ell)\sqrt{kpq}+\frac{2k}{\sqrt{kpq}}-2\ln(pq)}\\
&\le\E{Y_k}\cdot\pr{\frac{e^{2+2\sqrt{kpq}}kn(kpq+1)}{(k-\ell)^2}q^{-\frac{k+\ell-1}{2}}}^{k-\ell}
e^{\frac{2k}{\sqrt{kpq}}+O(\ln{k})}\\
&\le\E{Y_k}\cdot\pr{\frac{e^{O(\sqrt{\ln{np}})}kn\ln{np}}{q^{k}(k-\ell)^2}q^{\frac{k-\ell}{2}}}^{k-\ell}
e^{\frac{2k}{\sqrt{kpq}}+O(\ln{k})}\\
&\le\E{Y_k}\cdot\pr{\frac{kn}{(np)^{1.99}(k-\ell)^2}q^{\frac{k-\ell}{2}}}^{k-\ell}
e^{\frac{2k}{\sqrt{kpq}}+O(\ln{k})}\\
&=\E{Y_k}\cdot{g(k-\ell,kn/(np)^{1.99},q^{1/2},2)}
e^{\frac{2k}{\sqrt{kpq}}+O(\ln{k})},
\end{aligned}
$$
where in the end we used $g(x,y,z,\alpha)=\pr{\frac{yz^x}{x^{\alpha}}}^x$ from Lemma \ref{Lemma_g}.

Since $0<k-\ell<k/\sqrt{kpq}+1$, using \eqref{eq:function_g_upper_bound_0}, we get
$$
\begin{aligned}
g(k/\sqrt{kpq}&+1,kn/(np)^{1.99},q^{1/2},2)
\le\pr{\frac{O(np)}{(np)^{1.99}}q^{O\pr{\frac{k}{\sqrt{kpq}}}}}^{\frac{k}{\sqrt{kpq}}+1}
<\pr{\frac{1}{(np)^{0.98}}}^{\frac{k}{\sqrt{kpq}}+1}<1,\\
g(k-\ell,kn&/(np)^{1.99},q^{1/2},2)
\le\max\br{1,\exp\pr{\sqrt{\frac{kn}{(np)^{1.99}}}}}
<e^{\frac{k}{\sqrt{kpq}}}.
\end{aligned}
$$
Thus, by using Lemma \ref{Lemma_g},
$$
\frac{F_{\ell}}{\pr{\E{Y_k}}^2}\le
\frac{e^{\frac{3k}{\sqrt{kpq}}+O(\ln{k})}}{\pr{\frac{1-\varepsilon/2}{1-\varepsilon}}^{k}}\le
\pr{1-\varepsilon/2}^{k},
$$
$$
\frac{\D{Y_k}}{(\E{Y_k})^2}\le
\sum\limits_{\ell=1}^{k-\ceil{\frac{k}{\sqrt{kpq}}-\frac{1}{pq}}}\frac{F_{\ell}}{\pr{\E{Y_k}}^2}+o(1).
$$

2. In this case $\ell\le{k}\le\pr{k-\ell+\frac{1}{pq}}^{2}pq$. We apply the last bound below:

\begin{multline}\label{bound_sum_case2}
\sum\limits_{m=0}^{\ell-1}\frac{1}{((m+1)!)^2}
\pr{k-\ell+\frac{1}{pq}}^{2m}
(\ell{pq})^{m}
\le\frac{l^{l-1}}{l!}\sum\limits_{m=0}^{\ell-1}\frac{1}{(m+1)!}
\pr{k-\ell+\frac{1}{pq}}^{2m}
({pq})^{m}\\
\le\frac{e^{\ell}}{\ell}\sum\limits_{m=0}^{\ell-1}\frac{1}{(m+1)!}
\pr{k-\ell+\frac{1}{pq}}^{2m}(pq)^{m}
\le\frac{e^{\ell}}{\ell}\cdot\frac{\ell}{\ell!}\pr{k-\ell+\frac{1}{pq}}^{2\ell-2}(pq)^{\ell-1}.
\end{multline}

\begin{multline*}
F_{\ell}=\binom{n}{\ell,k-\ell,k-\ell,n-2k+\ell}
\pr{k+\frac{1}{pq}}^{2(k-\ell)}\ell^{\ell}{(pq)^{2k-\ell-2}}
q^{\binom{\ell}{2}-2\binom{k}{2}}\times\\
\times\sum\limits_{m=0}^{\ell-1}\frac{(\ell{pq})^{m+1}}{((m+1)!)^2}
\pr{k-\ell+\frac{1}{pq}}^{2m}\\
=\binom{n}{\ell,k-\ell,k-\ell,n-2k+\ell}
\pr{k+\frac{1}{pq}}^{2(k-\ell)}\ell^{\ell+1}{(pq)^{2k-\ell-1}}
q^{\binom{\ell}{2}-2\binom{k}{2}}\times\\
\times\sum\limits_{m=0}^{\ell-1}\frac{(\ell{pq})^{m}}{((m+1)!)^2}
\pr{k-\ell+\frac{1}{pq}}^{2m}\le
\end{multline*}

Use bound (\ref{bound_sum_case2}),
$$
\begin{aligned}
&\le\binom{n}{\ell,k-\ell,k-\ell,n-2k+\ell}
\pr{k+\frac{1}{pq}}^{2(k-\ell)}\ell^{\ell+1}{(pq)^{2k-\ell-1}}
q^{\binom{\ell}{2}-2\binom{k}{2}}\frac{e^{\ell}}{\ell !}\pr{k-\ell+\frac{1}{pq}}^{2\ell-2}(pq)^{\ell-1}\\
&=\binom{n}{\ell,k-\ell,k-\ell,n-2k+\ell}
\pr{k+\frac{1}{pq}}^{2(k-\ell)}(pq)^{2k-2}
q^{\binom{\ell}{2}-2\binom{k}{2}}
\frac{(e\ell)^{\ell}}{(\ell-1)!}
\pr{k-\ell+\frac{1}{pq}}^{2\ell-2}\\
&=\binom{n}{\ell,k-\ell,k-\ell,n-2k+\ell}
\pr{k+\frac{1}{pq}}^{2k-2}(pq)^{2k-2}
q^{\binom{\ell}{2}-2\binom{k}{2}}
\frac{(e\ell)^{\ell}}{(\ell-1)!}
\pr{1-\frac{\ell}{k+\frac{1}{pq}}}^{2\ell-2}\\
&\le\binom{n}{\ell,k-\ell,k-\ell,n-2k+\ell}
(kpq+1)^{2k-2}
q^{\binom{\ell}{2}-2\binom{k}{2}}
e^{2\ell}\ell
\pr{1-\frac{\ell}{k+\frac{1}{pq}}}^{2\ell-2}\\
&=\pr{\E{Y_k}}^2\cdot\frac{\binom{k}{\ell}\binom{n-k}{k-\ell}}{\binom{n}{k}}q^{\binom{\ell}{2}}
e^{2\ell}\ell\pr{1-\frac{\ell}{k+\frac{1}{pq}}}^{2\ell-2}.
\end{aligned}
$$
Using \eqref{eq:function_g_upper_bound_0}, we infer

\begin{align}\label{Fl/(EY)^2}
&\frac{F_{\ell}}{\pr{\E{Y_k}}^2}
\le\frac{\binom{k}{\ell}\binom{n-k}{k-\ell}}{\binom{n}{k}}
q^{\binom{\ell}{2}}
e^{2\ell}\ell.
\end{align}

Hence,

\begin{multline*}
\frac{F_{\ell}}{\pr{\E{X_k}}^2}\le \binom{k}{\ell}\frac{k!}{(k-\ell)!}
\frac{\frac{(n-k)!}{(n-2k+\ell)!}}{\frac{n!}{(n-k)!}}
q^{\frac{k}{2}\frac{\ell^2}{k}}
\ell e^{2\ell}
\le\pr{\frac{ek^2}{\ell}}^{\ell}
\frac{\frac{(n-k)!}{(n-2k+\ell)!}}{\frac{n!}{(n-k)!}}
q^{\frac{k}{2}\frac{\ell^2}{k}}
\ell e^{2\ell}\\
\le\pr{\frac{ek^2}{\ell}}^{\ell}
\frac{n^{k-\ell}}{n^{k}\pr{1-\frac{k}{n}}^k}
(O(np))^{\frac{\ell^2}{k}}
\ell e^{2\ell}
\le\pr{O\pr{\frac{k}{n}}\frac{(np)^{\frac{\ell}{k}}}{\frac{\ell}{k}}}^{\ell}
e^{\frac{k^2}{n-k}}\\
=\pr{g\pr{\frac{\ell}{k},O\pr{\frac{k}{n}},np,1}}^{k}
e^{O\pr{1+\frac{k^2}{n}}},
\end{multline*}

where in the end we used $g(x,y,z,\alpha)=\pr{\frac{yz^x}{x^{\alpha}}}^x$ from Lemma \ref{Lemma_g}. So, due to  Lemma \ref{Lemma_g},

\begin{multline}\label{Fl/(EY)^2 v2}
\frac{F_{\ell}}{\pr{\E{X_k}}^2}\le\pr{\max\br{e^{O\pr{\frac{k}{n}}},g\pr{1-\frac{1}{2\sqrt{kpq}},O\pr{\frac{k}{n}},np,1}}}^{k}
e^{O\pr{1+\frac{k^2}{n}}}=e^{O\pr{1+\frac{k^2}{n}}},
\end{multline}

since

$$
\begin{aligned}
&g\pr{1-\frac{1}{2\sqrt{kpq}},O\pr{\frac{k}{n}},np,1}\\
&=\pr{O\pr{\frac{\ln{np}}{(np)^{\frac{1}{2\sqrt{kpq}}}}}}^{1-\frac{1}{2\sqrt{kpq}}}\\
&=\exp\pr{\pr{{1-\frac{1}{2\sqrt{kpq}}}}\pr{\ln\ln{np}-e^{O(1)}\sqrt{\ln{np}}}}\\
&=\exp\pr{-e^{O(1)}\sqrt{\ln{np}}}.
\end{aligned}
$$

Next, let's consider two cases: $p\ge(\ln{n})^2/\sqrt{n}$ and $p<(\ln{n})^2/\sqrt{n}$. 

If $p\ge(\ln{n})^2/\sqrt{n}$, then $k^2=o(n)$.
Let $\ell_0:=\floor{\ln{n}/\ln\ln\ln{n}}$. Then, for $\ell\le\ell_0$, (\ref{Fl/(EY)^2}) does not exceed
$$
\begin{aligned}
\frac{F_{\ell}}{\pr{\E{X_k}}^2}
&\le\frac{\binom{k}{\ell}\binom{n-k}{k-\ell}}{\binom{n}{k}}
q^{\binom{\ell}{2}}
e^{2\ell}\ell\le\\
&\pr{\frac{k!}{(k-\ell)!}}^2
\frac{\frac{(n-k)!}{(n-2k+\ell)!}}{\frac{n!}{(n-k)!}}
\exp\pr{\frac{\ell^2}{2}\ln q+2\ell+\ln\ell+O(1)}\le\\
&\pr{\frac{k!}{(k-\ell)!}}^2
\frac{\frac{(n-k)!}{(n-2k+\ell)!}}{\frac{n!}{(n-k)!}}
\exp\pr{\frac{\ell^2}{2}pq+2\ell+\ln\ell+O(1)}\le\\
&
k^{2\ell}\frac{n^{k-\ell}}{n^k}
e^{\ell(p\ln{n}+O(1))}
\le\pr{\frac{O(k^2)}{n^{1-p}}}^{\ell}
\le\pr{\frac{O((\ln{n})^2)}{n^{1-p}p^2}}^{\ell}.
\end{aligned}
$$
Since, as can be easily verified, the expression $\frac{1}{n^{1-p}p^2}$ is convex w.r.t. the variable $p$,
$$
\begin{aligned}
\frac{F_{\ell}}{\pr{\E{X_k}}^2}
\le\pr{O\pr{\max\br{\frac{\ln^2 n}{n^{-\frac{(\ln{n})^2}{\sqrt{n}}}(\ln{n})^4},\frac{(\ln{n})^2}{n^{\varepsilon}(1-\varepsilon)^2}}}}^{\ell}
\le\pr{\frac{O(1)}{(\ln{n})^2}}^{\ell}
<\pr{\frac{1}{\ln{n}}}^{\ell}.
\end{aligned}
$$
For $\ell>\ell_0$, using \eqref{eq:function_g_upper_bound_1}, we get
$$
\begin{aligned}
&\pr{g\pr{\frac{\ell}{k},O\pr{\frac{k}{n}},np,1}}^{k}\\
&\le\pr{\max\br{\pr{O\pr{\frac{k^2}{n}}}^{\frac{\ell_0}{2k}},g\pr{\frac{\ell_0}{k},O\pr{\frac{k}{n}},np,1},g\pr{1-\frac{1}{2\sqrt{kpq}},O\pr{\frac{k}{n}},np,1}}}^k\\
&\le\max\br{\pr{O\pr{\frac{k^2}{n}(np)^{\frac{\ell_0}{k}}}}^{\ell_0},\exp\pr{-e^{O(1)}k\sqrt{\ln{n}}}}\\
&\le\max\br{\pr{\frac{O((\ln{np})^2)}{n^{1-p}p^2}}^{\ell_0},\exp\pr{-e^{O(1)}k\sqrt{\ln{n}}}}\\
&\le\max\br{\pr{\frac{1}{\ln{n}}}^{\frac{\ln{n}}{\ln\ln\ln{n}}},\exp\pr{-e^{O(1)}k\sqrt{\ln{n}}}}
=o\pr{\frac{1}{k}},
\end{aligned}
$$
So, bound (\ref{Fl/(EY)^2 v2}) does not exceed
$$
\frac{F_{\ell}}{\pr{\E{X_k}}^2}
\le\pr{g\pr{\frac{\ell}{k},O\pr{\frac{k}{n}},np,1}}^{k}
e^{O\pr{1+\frac{k^2}{n}}}
=o\pr{\frac{1}{k}}.
$$
Thus,
$$
\frac{\D{X_k}}{(\E{X_k})^2}\le
\sum\limits_{\ell=1}^{k}\frac{F_{\ell}}{\pr{\E{Y_k}}^2}\le
\sum\limits_{\ell=1}^{\ell_0}
\pr{\frac{1}{\ln{n}}}^{\ell}+
k\frac{\max\limits_{\ell>\ell_0}F_{\ell}}{\pr{\E{Y_k}}^2}
=o(1).
$$

If $p<(\ln{n})^2/\sqrt{n}$, then
$$
\frac{\ln{k}}{k(\ln{np})^4/np}\le
\frac{np(p\ln\ln{np}+p\ln(3/p))}{\varepsilon(\ln{np})^5}\le
\frac{(\ln{n})^4(\ln\ln{n}+\ln{n})}{\varepsilon(\frac{1}{2}\ln{n}-2\ln\ln{n})^5}=O(1),
$$
$$
\frac{k^2}{n}\le\frac{3\ln{np}}{np}k,
$$
$$
\frac{\D{X_k}}{(\E{X_k})^2}\le
\sum\limits_{\ell=1}^{k}\frac{F_{\ell}}{\pr{\E{X_k}}^2}\le
\exp\pr{O\pr{\frac{k^2}{n}}+O\pr{\ln{k}}}\le
\exp\pr{O\pr{\frac{(\ln{np})^4}{np}k}}.
$$

\section{Proof of Theorem \ref{th:concentration_max_ind_rooted_forests}}
The upper bound is given by Lemma~\ref{lemma:expectation_n_ind_forests_at_k_pm_eps} and the first moment method.
If $p\ge(\ln{n})^2/\sqrt{n}$, then the lower bound follows from Lemma~\ref{lemma:variance_n_ind_forests_at_k_m_eps} applied for $k=k_{-\varepsilon}$ and the second moment method.

If $p<(\ln{n})^2/\sqrt{n}$, then we will obtain the lower bound using Talagrand's inequality from \cite{Janson2000}. For a graph $G$ on the vertex set $[n]$ with the set of edges $\bigsqcup\limits_{i=1}^{n}x_i$ defined by the variables $x_i\subset\setdef{\br{i,j}}{j\in\irange{1}{i-1}}$, $i\in\irange{2}{n}$,
let $\varphi(G)=\varphi(x_2,\ldots,x_n)$ be the size of the largest induced forest in $G$.
First, notice that $\varphi(x_2,\ldots,x_n)$ satisfies Lipschitz condition. Indeed, $\varphi(x_2,\ldots,x_n)$ cannot decrease if an edge is removed. On the other hand, if any number of edges incident to a single vertex is added, $\varphi(x_2,\ldots,x_n)$ can decrease by at most $1$. Secondly, notice that $\varphi(x_2,\ldots,x_n)$ is certifiable by the function $\psi(r)=\ceil{r}$. Indeed, for any integer $r$, there exists a certificate of size $r$ enforcing the inequality $\varphi(x_2,\ldots,x_n)\ge{r}$, which is simply the set of vertices of the corresponding largest induced forest. Thus, Talagrand's inequality is applicable to the random variable $F(G(n,p))=\varphi(G(n,p))$:
$$
\begin{aligned}
\Pb{F(G(n,p))\le{k_{-\varepsilon}-1}}\Pb{F(G(n,p))\ge{k_{-\varepsilon/2}}}
&\le\exp\pr{-\frac{(k_{-\varepsilon/2}-k_{-\varepsilon}+1)^2}{4k_{-\varepsilon/2}}}\\
&\le\exp\pr{-\frac{\varepsilon^2}{4k_{-\varepsilon/2}p^2}}\\
&\le\exp\pr{-\frac{\varepsilon^{2}k_{-\varepsilon/2}}{40(\ln{np})^2}}.
\end{aligned}
$$
Using Lemma~\ref{lemma:variance_n_ind_forests_at_k_m_eps} for $k=k_{-\varepsilon/2}$ and Paley--Zygmund inequality, we get
$$
\Pb{F(G(n,p))\ge{k_{-\varepsilon/2}}}=\Pb{X_{k_{-\varepsilon/2}}>0}\ge\frac{(\op{E}{X_{k_{-\varepsilon/2}}})^2}{\op{E}{X_{k_{-\varepsilon/2}}^2}}
\ge\exp\pr{-\frac{k_{-\varepsilon/2}}{(\ln{np})^5}}
$$
and, finally, recalling that $p<(\ln{n})^2/\sqrt{n}$,
$$
\Pb{F(G(n,p))\le{k_{-\varepsilon}-1}}
\le\exp\pr{\frac{k_{-\varepsilon/2}}{(\ln{np})^5}-\frac{\varepsilon^{2}k_{-\varepsilon/2}}{40(\ln{np})^2}}
<\exp\pr{-\frac{\sqrt{n}}{(\ln{n})^4}}
\to0.
$$
\qed{}

\subsection*{Acknowledgements}

The research work of Kozhevnikov is supported by Russian Science Foundation, project 21‐71‐ 10092.
The research work of Akhmejanova is supported by King Abdullah University of Science and Technology (KAUST). The first author is a Young Russian Mathematics award winner (2020) and would like to thank its sponsors and jury.


\renewcommand{\refname}{References}

\section{Appendix}

\subsection{Coding trees with a fixed independent set}
\label{sub_trees_with_independent set}
It is known that a labelled tree on vertices $\irange{1}{n}$ can be encoded by a Pr\"{u}fer sequence $\pr{a_1,\ldots,a_{n-2}}\in\irange{1}{n}^{n-2}$. Let's fix some $m\in\irange{1}{n}$. Then the standard encoding procedure from \cite{Moon1970} can be slightly modified to encode only the trees with the independent set $\irange{1}{m}$ by two sequences $a=\pr{a_1,\ldots,a_{m-1}}\in\irange{m+1}{n}^{m-1}$ and $b=\pr{b_1,\ldots,b_{n-m-1}}\in\irange{1}{n}^{n-m-1}$.
Take any such tree and successively prune it's leaves, starting from leaves with smallest labels, as in the standard encoding procedure, until only one edge remains. If, however, the pruned leaf has the label in $\irange{1}{m}$, then write the label of its neighbour to the sequence $a$, otherwise, to the sequence $b$.

This modified encoding allows to find the number of such labelled trees on vertices $\irange{1}{n}$ with a given degree sequence $\pr{d_1,\ldots,d_n}$, which is exactly
\begin{equation}
\label{eq_n_trees_with_independent_set_given_deg}
\binom{n-m-1}{d_1-1,d_2-1,\ldots,d_m-1,n-1-\pr{d_1+\ldots,d_m}}
\binom{n+m-2-\pr{d_1+\ldots,d_m}}{d_{m+1}-1,\ldots,d_n-1}.
\end{equation}

From \eqref{eq_n_trees_with_independent_set_given_deg} it is evident that such a tree with the degree sequence $\pr{d_1,\ldots,d_n}$ exists iff $d_1+\ldots+d_n=2n-2$ and $d_1+\ldots+d_m\le{n-1}$.

\subsection{Counting forests with a fixed induced subforest}
\label{sub_forests_with_indiced_subforest}

Let $K_n$ be the complete graph on $n$ vertices. For any $F\subset{K_n}$, let's define $f(n,h,F)$ as the number of  rooted (with a distinguished vertex in each component) spanning forests of $K_n$ with $h$ components containing $F$ as an induced subgraph.

Let's consider a forest $F\subset{K_n}$ with $m$ components with numbers of vertices $f_1,\ldots,f_m$ and let $l:=\abs{V(F)}=f_1+\ldots+f_m$. Each $h$-component rooted spanning forest of $K_n$ containing $F$ as an induced subgraph can be uniquely constructed as follows:
\begin{itemize}
    \item contract each component of $F$ into a single vertex;
    \item label the contracted components by $\irange{1}{m}$ and the rest of the vertices of $K_{n+1}$ by $\irange{m+1}{n-\ell+m}$;
    \item add one more vertex labelled $n-\ell+m+1$;
    \item construct a labelled tree on the vertices $\irange{1}{n-\ell+m+1}$ in which $\deg(n-\ell+m+1)=h$;
    \item expand the vertices $\irange{1}{m}$ (which are the contracted components);
    \item remove the vertex labelled $n-\ell+m+1$.
\end{itemize}
Thus, applying \eqref{eq_n_trees_with_independent_set_given_deg},
\begin{multline*}
f(n,h,F)=\sum\limits_{d_1,\ldots,d_{n-\ell+m}=1}^{+\infty}
\binom{n-\ell}{d_1-1,d_2-1,\ldots,d_m-1,n-\ell+m-\pr{d_1+\ldots,d_m}}\cdot\\
\cdot\binom{n-\ell+2m-1-\pr{d_1+\ldots,d_m}}{d_{m+1}-1,\ldots,d_{n-\ell+m}-1,h-1}\cdot{f_1^{d_1}\cdot\ldots\cdot{f_m^{d_m}}}\\
=\sum\limits_{k_0,k_1,\ldots,k_{m}=0}^{+\infty}
\binom{n-\ell}{k_0,k_1,\ldots,k_m}\cdot{f_1^{k_1+1}\cdot\ldots\cdot{f_m^{k_m+1}}}
\cdot\sum\limits_{k_{m+1},\ldots,k_{n-\ell+m}=0}^{+\infty}
\binom{k_0+m-1}{k_{m+1},\ldots,k_{n-\ell+m},h-1}\\
=f_1\cdot\ldots\cdot{f_m}\cdot
\sum\limits_{k_0=0}^{n-\ell}
\binom{n-\ell}{k_0}l^{n-\ell-k_0}
\cdot\binom{k_{0}+m-1}{h-1}(n-\ell)^{k_{0}+m-h}
\end{multline*}
It is easy to see that $f(n,h,F)$ depends only on the sizes of the components of $F$. Therefore, we can define
$$
f(k,h,\br{f_1,\ldots,f_m}):=f(k,h,F).
$$
Thus,
\begin{equation}
\label{eq_induced_forest_extension_bound}
f(k,h,\br{f_1,\ldots,f_m})=
f_1\cdot\ldots\cdot{f_m}
\sum\limits_{k_0=0}^{k-\ell}
\binom{k-\ell}{k_0}\ell^{k-\ell-k_0}
\cdot\binom{k_{0}+m-1}{h-1}(k-\ell)^{k_{0}+m-h}.
\end{equation}

\end{document}